\documentclass{amsart}
\usepackage{amsrefs}
\usepackage[all]{xy}
\usepackage[mathscr]{eucal}
\usepackage{graphicx,amssymb}
\usepackage{pinlabel}
\usepackage{placeins}
\usepackage{url}
\CompileMatrices

\newtheorem{theorem}{Theorem}

\theoremstyle{definition}
\newtheorem{definition}[theorem]{Definition}

\newtheorem*{examples}{Examples}

\theoremstyle{remark}

\numberwithin{equation}{section}

\begin{document}

\newcommand{\Diff}{\operatorname{Diff}}
\newcommand{\Homeo}{\operatorname{Homeo}}
\newcommand{\Hom}{\operatorname{Hom}}
\newcommand{\Exp}{\operatorname{Exp}}
\newcommand{\Orb}{\operatorname{\textup{Orb}}}
\newcommand{\Emb}{\operatorname{Emb}}
\newcommand{\Stab}{\operatorname{Stab}}

\newcommand{\COrb}{\operatorname{\star\textup{Orb}}}
\newcommand{\CROrb}{\operatorname{\scriptscriptstyle{\blacklozenge}\scriptstyle\textup{Orb}}}
\newcommand{\ssslozenge}{\scriptscriptstyle{\blacklozenge}}
\newcommand{\ssstriangledown}{\scriptscriptstyle{\blacktriangledown}}

\newcommand{\Stwo}{\mbox{$\displaystyle S^2$}}
\newcommand{\Sn}{\mbox{$\displaystyle S^n$}}
\newcommand{\supp}{\operatorname{supp}}
\newcommand{\intr}{\operatorname{int}}
\newcommand{\kernel}{\operatorname{ker}} 

\newcommand{\RR}{\mathbb{R}} 
\newcommand{\ZZ}{\mathbb{Z}}
\newcommand{\Beta}{\mathrm{B}}

\newcommand{\orbify}[1]{\ensuremath{\mathcal{#1}}}
\newcommand{\starfunc}[1]{\ensuremath{{}_\star{#1}}}
\newcommand{\lozengefunc}[1]{\ensuremath{{}_{\scriptscriptstyle{\blacklozenge}}{#1}}}
\newcommand{\redfunc}[1]{\ensuremath{{}_\bullet{#1}}}

\newcommand{\OrbDiff}{\ensuremath{\Diff_{\Orb}}}
\newcommand{\RedOrbDiff}{\ensuremath{\Diff_{\textup{red}}}}
\newcommand{\COrbDiff}{\ensuremath{\Diff_{\COrb}}}
\newcommand{\CROrbDiff}{\ensuremath{\Diff_{\CROrb}}}
\newcommand{\OrbMaps}{\ensuremath{C_{\Orb}}}
\newcommand{\RedOrbMaps}{\ensuremath{C_{\textup{red}}}}
\newcommand{\COrbMaps}{\ensuremath{C_{\COrb}}}
\newcommand{\CROrbMaps}{\ensuremath{C_{\CROrb}}}
\newcommand{\Frechet}{Fr\'{e}chet\ }
\newcommand{\Frechetnospace}{Fr\'{e}chet}
\newcommand*{\longhookrightarrow}{\ensuremath{\lhook\joinrel\relbar\joinrel\rightarrow}}

\title[On the Inheritance of Orbifold Substructures]{On the Inheritance of Orbifold Substructures}

\author{Joseph E. Borzellino}

\address{Department of Mathematics, California Polytechnic State
  University, 1 Grand Avenue, San Luis Obispo, California 93407}
\email{jborzell@calpoly.edu}
% \thanks will become a 1st page footnote.  \thanks{The first author
%   was supported in part by NSF Grant \#000000.}

% Information for second author
\author{Victor Brunsden} \address{Department of Mathematics and
  Statistics, Penn State Altoona, 3000 Ivyside Park, Altoona,
  Pennsylvania 16601} \email{vwb2@psu.edu}
% \thanks{Support information for the second author.}

% General info

\date{\today} \commby{Editor} 
\begin{abstract}
In a previous article, we defined a very flexible notion of suborbifold and characterized those suborbifolds which can arise as the images of orbifold embeddings. In particular, suborbifolds are images of orbifold embeddings precisely when they are \emph{saturated} and \emph{split}. This article addresses the problem of orbifold structure inheritance for three orbifolds $\orbify{Q}\subset\orbify{P}\subset\orbify{O}$. We identify an appealing but ultimately inadequate notion of an inherited canonical orbifold substructure. In particular, we give a concrete example where the orbifold structure of $\orbify{Q}$ is canonically inherited from $\orbify{P}$, and the orbifold structure of $\orbify{P}$ is canonically inherited from $\orbify{O}$, but the orbifold structure of $\orbify{Q}$ is not canonically inherited from $\orbify{O}$. On the other hand, it is easy to see that when $\orbify{Q}$ is embedded in $\orbify{P}$, and $\orbify{P}$ is embedded in $\orbify{O}$, all of the canonical inherited orbifold substructures will agree. We also investigate the property of saturation in this context, and give an example of a suborbifold with the canonical orbifold substructure that is not saturated.
\end{abstract}

\maketitle

\section{Introduction}\label{IntroSection}

In \cite{MR3426692}, we defined a very flexible notion of suborbifold and characterized those suborbifolds which can arise as the images of orbifold embeddings. In particular, suborbifolds are images of orbifold embeddings if and only if they are \emph{saturated} and \emph{split}. For manifolds, it is a fundamental result of differential topology that submanifolds are precisely the images of embeddings, and in fact, many authors use this characterization as the definition of submanifold. Thus, we were surprised to find examples of suborbifolds that were not images of orbifold embeddings. Because of these examples, we began to look at the issue of orbifold substructure inheritance for three orbifolds $\orbify{Q}\subset\orbify{P}\subset\orbify{O}$. To do this, we first look more deeply at the definition of suborbifold and identify an appealing but ultimately inadequate notion of an inherited canonical orbifold substructure. The appeal comes from the fact that certain choices in the definition of suborbifold can be made in a unique way. The inadequacy comes from the observation that these unique canonical suborbifold structures do not persist through inclusion. In particular, we give a concrete example of three orbifolds $\orbify{Q}\subset\orbify{P}\subset\orbify{O}$, where the orbifold structure of $\orbify{Q}$ is canonically inherited from $\orbify{P}$, and the orbifold structure of $\orbify{P}$ is canonically inherited from $\orbify{O}$, but the orbifold structure of $\orbify{Q}$ is not canonically inherited from $\orbify{O}$. On the other hand, it is easy to see that when $\orbify{Q}$ is embedded in $\orbify{P}$, and $\orbify{P}$ is embedded in $\orbify{O}$, that all of the canonical inherited orbifold substructures will agree. To be clear, by embedded, we mean realized as the image of an orbifold embedding \cite{MR3426692}. In some ways, this phenomenon can be considered an orbifold structure analog of the difference between one-to-one immersions and embeddings in differential topology. That is, the topology induced by the immersion may not agree with the relative topology of its image. In the orbifold context, one could say we  have ``immersed" orbifold substructures that are not ``embedded," although we do not use or define these terms. We also investigate further the property of saturation, and give an example of a suborbifold with the canonical orbifold substructure that is not saturated.

\section{Orbifold Preliminaries}
All of the following definitions come from our article \cite{MR3426692} and the references therein (especially, \cites{MR2523149,MR3117354}). This will be our standard reference for orbifold background material. We also assume the reader is familiar with the definition of a smooth orbifold modeled after Thurston \cite{Thurston78}. Such orbifolds are referred to as classical effective orbifolds in \cite{MR2359514}. As such, each point $x$ in a smooth orbifold $\orbify{O}$, has a neighborhood $U_x$ or orbifold chart $(\tilde U_x,\Gamma_x)$ or $(\tilde U_x,\Gamma_x, \rho_x, \phi_x)$ where $\tilde U_x\cong\RR^n$, $\phi_x:\tilde U_x/\Gamma_x\to U_x$ is a homeomorphism, and $\phi_x(0)=x$. In the 4-tuple  notation, we are making explicit the (faithful) representation $\rho_x:\Gamma_x\to\text{Diff}^\infty(\tilde U_x,0)$. $\text{Diff}^\infty(\tilde U_x,0)$ are the group of smooth diffeomorphisms that leave the origin fixed. The \emph{isotropy group of $x$} is the group $\Gamma_x$. It is unique up to isomorphism. By the Bochner-Cartan theorem \cites{MR0073104,MR2523149}, the smooth action of $\Gamma_{x}$ is smoothly conjugate to the linear action on $\tilde U_x$ given by the differential of the action. So without loss of generality, we may assume $\rho_x:\Gamma_x\to O(n)$.
These charts are subject to overlap compatibility conditions that give $\orbify{O}$ its orbifold structure. More detail can be found in \cite{MR2523149}. We now recall several definitions related to the notion of suborbifold from \cites{MR2973378,MR3426692}:

\begin{definition}\label{SubOrbifold}
  A  \emph{suborbifold} \orbify{P} of an orbifold \orbify{O} consists
  of the following:
  \begin{enumerate}
  \item A subspace $X_{\orbify{P}}\subset X_{\orbify{O}}$ equipped
    with the subspace topology.
  \item\label{FlexibleDef} For each $x\in X_{\orbify{P}}$ and neighborhood $W$ of $x$ in
    $X_{\orbify{O}}$ there is an orbifold chart $(\tilde U_x,
    \Gamma_x, \rho_x, \phi_x)$ about $x$ in \orbify{O} with
    $U_x\subset W$, a subgroup $\Lambda_x \subset \Gamma_x$ of the
    isotropy group of $x$ in \orbify{O} and a $\rho_x(\Lambda_x)$
    invariant submanifold $\tilde V_x\subset \tilde U_x \cong \RR^n$,
    so that $(\tilde V_x, \Lambda_x/\Omega_x,\rho_x\lvert_{\Lambda_x},\psi_x)$ is
    an orbifold chart for $\orbify{P}$, where $\Omega_x=\left\{\gamma\in\Lambda_x\mid \rho_x(\gamma)\lvert_{\tilde{V}_x}=\text{Id}\right\}$. (In particular, the \emph{intrinsic} isotropy subgroup at $x\in\orbify{P}$ is  $\Lambda_x/\Omega_x$).
  \item
      For $x$ in $\orbify{P}$, $V_x = \psi_x(\tilde V_x/\rho_x(\Lambda_x))
      =U_x\cap X_{\orbify{P}}$
    is an orbifold chart.
  \end{enumerate}
\end{definition}

Implicit in this definition is the requirement that the invariant submanifolds $\tilde{V}_x$ be smooth, and that the collection of charts $\{(\tilde V_x, \Lambda_x/\Omega_x,\rho_x\lvert_{\Lambda_x},\psi_x)\}$ satisfy the overlap compatibility conditions of an orbifold, thus giving $\orbify{P}$ the structure of a smooth orbifold. As previously noted, condition (\ref{FlexibleDef}) is not very restrictive (see \cite{MR3426692}).

Motivated by Thurston's notion of suborbifold \cite{Thurston78}, we made the following definition:

\begin{definition}\label{FullSuborbifoldDef} $\orbify{P}\subset\orbify{O}$ is a \emph{full suborbifold} of $\orbify{O}$ if $\orbify{P}$ is a suborbifold with
$\Lambda_x=\Gamma_x$ for all $x\in\orbify{P}$.
\end{definition}

When necessary for clarity, we will use the notation $\Gamma_{x,\orbify{O}}$ to denote the intrinsic isotropy group of a point $x$ in an orbifold $\orbify{O}$, and use the subscript $\orbify{O}$ as well on needed subgroups of $\Gamma_{x,\orbify{O}}$. When the base point $x$ is clear, we may drop it as well. In the case of a suborbifold $\orbify{P}\subset\orbify{O}$ we always have the following exact sequence of groups

$$1\longrightarrow\Omega_{x,\orbify{O}}\longrightarrow\Lambda_{x,\orbify{O}}\subset\Gamma_{x,\orbify{O}}\longrightarrow\Gamma_{x,\orbify{P}}\longrightarrow 1$$
where $\Gamma_{x,\orbify{P}}$ denotes the intrinsic isotropy group of $\orbify{P}$ at $x$.

In characterizing those suborbifolds that are images of orbifold embeddings in \cite{MR3426692}, we identified the following two conditions.

\begin{definition}\label{SplitSuborbifoldDef} We say that $\orbify{P}\subset\orbify{O}$ is a \emph{split} suborbifold of $\orbify{O}$ if the exact sequence above is (right) split for all $x\in\orbify{P}$. That is, there is a group homomorphism $\sigma:\Gamma_{x,\orbify{P}}\to\Lambda_{x,\orbify{O}}$ such that the composition $q\circ\sigma=\text{Id}$, where $q:\Lambda_{x,\orbify{O}}\to\Gamma_{x,\orbify{P}}$ is the quotient homomorphism:

\begin{equation*}
\xymatrix{1\ar[r] & \Omega_{x,\orbify{O}}\ar[r] & \Lambda_{x,\orbify{O}}\ar[r]^-{q} & \Gamma_{x,\orbify{P}}\ar@/^.75pc/[l]^\sigma\ar[r] & 1.}
\end{equation*}
\end{definition}

Note that if $\orbify{P}\subset\orbify{O}$ is split, we have $\Lambda_{x,\orbify{O}}\cong\Omega_{x,\orbify{O}}\rtimes\Gamma_{x,\orbify{P}}$, a semidirect product, and in the case that the groups are abelian $\Lambda_{x,\orbify{O}}\cong\Omega_{x,\orbify{O}}\times\Gamma_{x,\orbify{P}}$, the direct product. Of course, if $\Omega_{x,\orbify{O}}$ or $\Gamma_{x,\orbify{P}}$ is trivial, then $\orbify{P}$ is split as well.

\begin{definition}\label{SaturatedSuborbifoldDef} We say that $\orbify{P}\subset\orbify{O}$ is a \emph{saturated} suborbifold of $\orbify{O}$ if for each $x\in\orbify{P}$ and $\tilde y\in\tilde V_x$, we have that
$\left(\Gamma_{x,\orbify{O}}\cdot\tilde y\right)\cap\tilde V_x=\Lambda_{x,\orbify{O}}\cdot\tilde y$.
\end{definition}

The main result of \cite{MR3426692} proved that a suborbifold $\orbify{P}\subset\orbify{O}$ was the image of an orbifold embedding if and only if $\orbify{P}$ was both saturated and split. In this paper, we need the following definition to identify those suborbifolds that have orbifold substructures which are inherited in essentially a unique way.

\begin{definition}\label{CanSubstructure} A suborbifold $\orbify{P}\subset\orbify{O}$ has the \emph{canonical orbifold substructure}, $\orbify{P}^\orbify{O}$, \emph{inherited from} $\orbify{O}$, if, in definition~\ref{SubOrbifold}, $\Lambda_x=\text{Stab}^{\Gamma_x}(\tilde V_x)$ for all $x\in X_\orbify{P}$. That is, for all $x\in X_\orbify{P}$, $\Lambda_x=\{\gamma\in\Gamma_x\mid \gamma\cdot\tilde V_x\subset\tilde V_x\}$, the entire stabilizer subgroup of $\tilde V_x$ in $\Gamma_x$.
\end{definition}

Another way to think about this is that if $\orbify{P}$ is covered at $x$ by a chart $\tilde{V}_x\subset\tilde U_x$, then the group $\Lambda_x$ is completely determined when $\orbify{P}$ carries the canonical orbifold substructure inherited from $\orbify{O}$. 

It is also clear that a saturated suborbifold $\orbify{P}\subset\orbify{O}$ has the canonical orbifold substructure inherited from $\orbify{O}$. Since embedded suborbifolds are both saturated and split \cite{MR3426692}, and the composition of orbifold embeddings is an orbifold embedding, it follows that in the case $\orbify{Q}$ is embedded in $\orbify{P}$, and $\orbify{P}$ is embedded in $\orbify{O}$, that all of the canonical inherited orbifold substructures will agree.

\begin{examples}\label{CanonicalVsSaturatedEx} To illustrate some of these definitions, we revisit the following examples which are all in \cite{MR3426692}. It is clear that any full suborbifold carries a canonical suborbifold structure. However, this condition is not necessary as illustrated by examples~10, 11 and 14 of \cite{MR3426692}, all of which carry canonical suborbifold structures, but are not full suborbifolds. On the other hand, each of these suborbifolds are saturated. Furthermore, the suborbifolds of example~13, which are arguably the most counterintuitive (\cite{Weilandt2016}) examples presented in \cite{MR3426692}, do not carry inherited canonical orbifold substructures. They differ from the other examples in that none of these suborbifolds are saturated.
\end{examples}

Lastly, we'd like to point out that a recent preprint of Weilandt \cite{Weilandt2016} expanded on our work in \cite{MR3426692} and includes a detailed discussion of the challenges surrounding the notion of suborbifold. For Weilandt, however, all suborbifolds are saturated although his notion of full suborbifold turns out to be the same as ours. As pointed out earlier, saturated suborbifolds already possess the canonical orbifold substructure. In section~\ref{CanVsSatSection}, we show that a suborbifold can possess the canonical orbifold substructure, but not be saturated.

\section{Inherited Canonical Orbifold Substructures}\label{CanonicalSubstructures}
In this section, we show how to use group algebras to construct orbifolds $\orbify{Q}\subset\orbify{P}\subset\orbify{O}$ and give concrete expressions for the inherited canonical orbifold substructures.

\subsection{Review of Group Algebras} For a finite group $\Gamma$, let $\RR[\Gamma]$ denote the corresponding group algebra. That is, $\RR[\Gamma]$ is the $|\Gamma|$-dimensional $\RR$-vector space consisting of all formal sums $\sum_{\gamma\in\Gamma}c_\gamma\gamma$, with $c_\gamma\in\RR$, and the obvious scalar multiplication and component-wise addition. That is, for $\alpha=\sum c_\gamma\gamma$ and $\beta=\sum d_\gamma\gamma$:
$$\alpha+\beta=\sum_{\gamma\in\Gamma} (c_\gamma+d_\gamma)\gamma,\quad\text{and }k\alpha=\sum_{\gamma\in\Gamma} (kc_\gamma)\gamma,\ k\in\RR.$$
Because $\Gamma$ is a group, $\RR[\Gamma]$ becomes an algebra via:

\begin{gather*}\left(\sum_{\gamma\in\Gamma} c_\gamma\gamma\right)\left(\sum_{\delta\in\Gamma} d_\delta\delta\right)=\sum_{\gamma,\delta\in\Gamma}(c_\gamma d_\delta)\gamma\delta=\sum_{\nu\in\Gamma} e_{\nu}\nu,\\
\intertext{where}
e_\nu=\sum_{\gamma\delta=\nu}c_\gamma d_\delta=\sum_{\gamma\in\Gamma}c_\gamma d_{\gamma^{-1}\nu}.
\end{gather*}

There is a natural (left) action of $\Gamma$ on $\RR[\Gamma]$ given by component-wise conjugation:
$$\delta\cdot\alpha=\delta\cdot\sum_{\gamma\in\Gamma}c_\gamma\gamma=\sum_{\gamma\in\Gamma}c_\gamma(\delta\gamma\delta^{-1}).$$
This action is effective precisely when $\Gamma$ has a trivial center, as $\delta\gamma\delta^{-1}=\gamma$ for all $\gamma\in\Gamma$ implies $\delta\in C(\Gamma)$, the center of $\Gamma$. In particular, the action is effective when $\Gamma$ is a simple, nonabelian group.

Now, let $\Beta\subset\Gamma$ be a subgroup and consider its group algebra $\RR[\Beta]$. Define the stabilizer of $\RR[\Beta]$ in $\Gamma$ to be
$$\Stab^\Gamma(\RR[\Beta])=\{\gamma\in\Gamma\mid\gamma\cdot\RR[\Beta]\subset\RR[\Beta]\}.$$
Let $\beta=\sum_{\eta\in\Beta}c_\eta\eta\in\RR[\Beta]$ be arbitrary. Then, $\gamma\in\Stab^\Gamma(\RR[\Beta])$ implies $\gamma\cdot\beta\in\RR[\Beta]$ which in turn implies that $\gamma\eta\gamma^{-1}\in\Beta$ for all $\eta\in\Beta$. This means that $\gamma\in N_\Gamma(\Beta)$, the normalizer of $B$ in $\Gamma$. Thus, $\Stab^\Gamma(\RR[\Beta])=N_\Gamma(\Beta)$. Similarly, let
$$\Gamma^{\RR[\Beta]}=\{\gamma\in\Gamma\mid\gamma\cdot\beta=\beta\text{ for all }\beta\in\Beta\}.$$
Let $\beta=\sum_{\eta\in\Beta}c_\eta\eta\in\RR[\Beta]$ again be arbitrary. Then, $\gamma\in\Gamma^{\RR[\Beta]}$ implies that $\gamma\eta\gamma^{-1}=\eta$ for all $\eta\in\Beta$. This means that $\gamma\in C_\Gamma(\Beta)$, the centralizer of $B$ in $\Gamma$. Thus, $\Gamma^{\RR[\Beta]}=C_\Gamma(\Beta)$. Note also that $C_\Gamma(\Beta)$ is a normal subgroup of $N_\Gamma(\Beta)$.

\subsection{Constructing Orbifolds from the Group Algebra}\label{GroupAlgebraConstructionSection} We start by choosing a finite centerless group $\Gamma$, and two additional nontrivial subgroups $\Delta$ and $\Beta$ each properly contained in the other so that $\{e\}\subsetneq\Delta\subsetneq\Beta\subsetneq\Gamma$. Let $\tilde U=\RR[\Gamma]$, $\tilde V=\RR[\Beta]$, and $\tilde W=\RR[\Delta]$. Naturally, $\tilde{W}\subsetneq\tilde{V}\subsetneq\tilde{U}$. By choice of $\Gamma$, $\orbify{O}=\tilde U/\Gamma$, is an $|\Gamma|$-dimensional (classical effective) orbifold covered by a single orbifold chart $(\tilde U,\Gamma)$, and the isotropy group of the origin is $\Gamma$. We now build a suborbifold $\orbify{P}$ of $\orbify{O}$ using definition~\ref{SubOrbifold}. We first choose our subgroup $\Lambda$ of $\Gamma$ to be $\Lambda=\Stab^\Gamma(\tilde V)=N_\Gamma(\Beta)$. In this case, it follows that $\Omega=C_\Gamma(\Beta)$. If we let $\Gamma_{\orbify{P}}=\Lambda/\Omega=N_\Gamma(\Beta)/C_\Gamma(\Beta)$, then $\orbify{P}=\tilde V/\Gamma_{\orbify{P}}$ is a suborbifold of $\orbify{O}$ equipped with the canonical orbifold substructure, $\orbify{P}^{\orbify{O}}$, inherited from $\orbify{O}$. The intrinsic isotropy group (at the origin) is equal to $\Gamma_{\orbify{P}}$.

The interesting observation is that now we have the option of equipping $\orbify{Q}$ with either of two canonical suborbifold structures; one inherited from $\orbify{O}$ and one inherited from $\orbify{P}$. Specifically, we can build a suborbifold $\orbify{Q}=\tilde W/\Gamma_{\orbify{Q}}^{\orbify{O}}$ of $\orbify{O}$, equipped with the canonical orbifold substructure, $\orbify{Q}^{\orbify{O}}$, inherited from $\orbify{O}$. In this case, $\Gamma_{\orbify{Q}}^{\orbify{O}}=N_\Gamma(\mathrm{\Delta})/C_\Gamma(\mathrm{\Delta})$ is the intrinsic isotropy group (at the origin) of $\orbify{Q}^{\orbify{O}}$. Or, we can build a suborbifold $\orbify{Q}=\tilde W/\Gamma_{\orbify{Q}}^{\orbify{P}}$ of $\orbify{P}$, equipped with the canonical orbifold substructure, $\orbify{Q}^{\orbify{P}}$, inherited from $\orbify{P}$. In this case, $\Gamma_{\orbify{Q}}^{\orbify{P}}=\Stab^{\Gamma_{\orbify{P}}}(\tilde W)/\Gamma_{\orbify{P}}^{\tilde W}$ is the intrinsic isotropy group (at the origin) of $\orbify{Q}^{\orbify{P}}$. Here, $\Gamma_{\orbify{P}}^{\tilde W}$ represents the subgroup of $\Gamma_{\orbify{P}}$ that fixes $\tilde W$ pointwise. One of our main goals is to construct an example of a suborbifold $\orbify{Q}\subset\orbify{P}\subset\orbify{O}$, where $\orbify{P}$ has the canonical orbifold structure inherited from $\orbify{O}$, and $\orbify{Q}$ has the canonical orbifold structure inherited from $\orbify{P}$, but $\orbify{Q}$ does not have the canonical orbifold structure inherited from $\orbify{O}$. In other words, we produce an example where $\Gamma_{\orbify{Q}}^{\orbify{O}}$ differs from $\Gamma_{\orbify{Q}}^{\orbify{P}}$, thus showing that canonical orbifold substructures are not uniquely inherited in general. Before we construct our example, we will unwind the definition of $\Gamma_{\orbify{Q}}^{\orbify{P}}$ for reference, even though we do not need its full generality in our construction.

\subsubsection{Unwinding $\Gamma_{\orbify{Q}}^{\orbify{P}}$}\label{UnwindingSubsection} To unwind $\Stab^{\Gamma_{\orbify{P}}}(\tilde W)/\Gamma_{\orbify{P}}^{\tilde W}$, it is best to think geometrically. First, we get an expression for $\Stab^{\Gamma_{\orbify{P}}}(\tilde W)$. $\Gamma_{\orbify{P}}$ consists of those elements of $\Gamma$ that leave $\tilde V$ invariant, and where we regard two such elements to be the same if they differ by an element of $\Gamma$ that leaves $\tilde V$ pointwise fixed. Since $\Stab^{\Gamma_{\orbify{P}}}(\tilde W)$ consists of the elements of $\Gamma_{\orbify{P}}$ that leave $\tilde W\subset\tilde V$ invariant, we see that such elements must be those elements of $\Gamma$ that leave both $\tilde W$ and $\tilde V$ invariant, and where we regard two elements the same if they differ by an element of $\Gamma$ that leaves $\tilde V$ pointwise fixed. That is,
$$\Stab^{\Gamma_{\orbify{P}}}(\tilde W)=\big[\Stab^{\Gamma}(\tilde W)\cap\Stab^{\Gamma}(\tilde V)\big]\big /\Gamma^{\tilde V}=\left[N_\Gamma(\Delta)\cap N_\Gamma(\Beta)\right]/C_\Gamma(\Beta).$$
Now $\Gamma_{\orbify{P}}^{\tilde W}$ consists of those elements of $\Gamma_{\orbify{P}}$ which leave $\tilde W$ pointwise fixed. Interpreting $\Gamma_{\orbify{P}}$ as before, we can conclude that
$$\Gamma_{\orbify{P}}^{\tilde W}=\big[\Stab^{\Gamma}(\tilde V)\cap\Gamma^{\tilde W}\big] \big /\Gamma^{\tilde V}=\left[N_\Gamma(\Beta)\cap C_\Gamma(\Delta)\right]/C_\Gamma(\Beta).$$
By the third isomorphism theorem we have:
\begin{align*}
\Gamma_{\orbify{Q}}^{\orbify{P}}&=\Stab^{\Gamma_{\orbify{P}}}(\tilde W)/\Gamma_{\orbify{P}}^{\tilde W}\\
&=\big[\Stab^{\Gamma}(\tilde W)\cap\Stab^{\Gamma}(\tilde V)\big]\big /\Gamma^{\tilde V} \bigg / \big[\Stab^{\Gamma}(\tilde V)\cap\Gamma^{\tilde W}\big] \big /\Gamma^{\tilde V}\\
&=\left[N_\Gamma(\Delta)\cap N_\Gamma(\Beta)\right]/C_\Gamma(\Beta)\bigg / \left[N_\Gamma(\Beta)\cap C_\Gamma(\Delta)\right]/C_\Gamma(\Beta)\\
&\cong \left[N_\Gamma(\Delta)\cap N_\Gamma(\Beta)\right] \big/ \left[N_\Gamma(\Beta)\cap C_\Gamma(\Delta)\right]
\end{align*}
Since $C_\Gamma(\Delta)\subset N_\Gamma(\Delta)$, the last expression may be written:
\begin{align*} \Gamma_{\orbify{Q}}^{\orbify{P}}&=
\left[N_\Gamma(\Delta)\cap N_\Gamma(\Beta)\right] \big/ \left[ (N_\Gamma(\Delta)\cap N_\Gamma(\Beta)) \cap C_\Gamma(\Delta)\right]\\
&\cong \left[(N_\Gamma(\Delta)\cap N_\Gamma(\Beta))\cdot C_\Gamma(\Delta)\right]\big / C_\Gamma(\Delta),
\end{align*}
where the last isomorphism is obtained by the second isomorphism theorem after noting that $C_\Gamma(\Delta)$ is a normal subgroup of
$N_\Gamma(\Delta)$, and thus a normal subgroup of $N_\Gamma(\Delta)\cap N_\Gamma(\Beta)$.

\section{Construction of Incompatible Inherited Canonical Orbifold Substructures}\label{MainCounterExample} We now return to the question of incompatibility of inherited canonical orbifold substructures. Since this question is of a local nature, it is sufficient to consider orbifolds that arise as quotients of a single orbifold chart. We will use the techniques of the previous section to construct three orbifolds, $\orbify{Q}\subset\orbify{P}\subset\orbify{O}$, where the canonical orbifold substructure that $\orbify{Q}$ inherits from $\orbify{P}$ as a suborbifold is different than the canonical orbifold substructure that $\orbify{Q}$ inherits from $\orbify{O}$ as a suborbifold. 

As mentioned at the end of paragraph preceding section~\ref{UnwindingSubsection}, we want to exhibit a situation where $\Gamma_{\orbify{Q}}^{\orbify{O}}$ differs from $\Gamma_{\orbify{Q}}^{\orbify{P}}$. Given the group algebra setup and computations in the previous sections, it is sufficient to find groups $\Delta\subset\Beta\subset\Gamma$, where
$$\Gamma_{\orbify{Q}}^{\orbify{P}}\cong\left[(N_\Gamma(\Delta)\cap N_\Gamma(\Beta))\cdot C_\Gamma(\Delta)\right]\big / C_\Gamma(\Delta)\ncong N_\Gamma(\Delta)/C_\Gamma(\Delta)\cong\Gamma_{\orbify{Q}}^{\orbify{O}}.$$

To that end, let $A_k$ denote the alternating group of degree $k$. For our groups, we let $\Gamma=A_5$, $\Beta=A_4$, and $\Delta=A_3\cong\ZZ_3$. Thus, $\tilde U=\RR[A_5]$, $\tilde V=\RR[A_4]$, and $\tilde W=\RR[\ZZ_3]$, vector spaces of dimensions 60, 12, and 3, respectively. We first compute $\Gamma_{\orbify{Q}}^{\orbify{O}}\cong N_\Gamma(\Delta)/C_\Gamma(\Delta)$. We use standard facts from group theory and the Groupprops Subwiki \cite{Groupprops} as our reference. Since the order of the conjugacy class of $\Delta$ in $\Gamma$ is 10, and equal to $[\Gamma:N_\Gamma(\Delta)]$, the index of the normalizer of $\Delta$ in $\Gamma$, we find that $N_\Gamma(\Delta)$ has order 6. The only subgroup of order 6 in $A_5$ is the twisted symmetric group of degree 3, $S_3\cong\ZZ_2\ltimes\ZZ_3$. The centralizer, $C_\Gamma(\Delta)$ of $\Delta$ in $\Gamma$ is a normal subgroup of $N_\Gamma(\Delta)$, and thus must be $\Delta$ since $\Delta$ is abelian and the centralizer cannot be all of $S_3$. We conclude that 
$$\Gamma_{\orbify{Q}}^{\orbify{O}}\cong N_\Gamma(\Delta)/C_\Gamma(\Delta)\cong(\ZZ_2\ltimes\ZZ_3)/\ZZ_3\cong\ZZ_2.$$
This means that in the canonical orbifold substructure that $\orbify{Q}$ inherits from $\orbify{O}$, the origin has non-trivial isotropy $\ZZ_2$.

Next, we compute $\Gamma_{\orbify{P}}^{\orbify{O}}\cong N_\Gamma(\Beta)/C_\Gamma(\Beta)$. Standard results from group theory show that the order of the conjugacy class of $\Beta$ in $\Gamma$ is 5. This implies that the order of the normalizer $N_\Gamma(\Beta)$ has order 12. Since $\Beta=A_4$ has order 12, we conclude that $N_\Gamma(\Beta)=\Beta$. Again, since the centralizer, $C_\Gamma(\Beta)$, of $\Beta$ in $\Gamma$ is a normal subgroup of $N_\Gamma(\Beta)=\Beta$, we find that $C_\Gamma(\Beta)=C(\Beta)$, the center of $\Beta$. Since $A_4$ is centerless, $C_\Gamma(\Beta)=\{1\}$, and 
thus,
$$\Gamma_{\orbify{P}}^{\orbify{O}}\cong N_\Gamma(\Beta)/C_\Gamma(\Beta)\cong\Beta/\{1\}\cong A_4.$$

Lastly, we compute $\Gamma_{\orbify{Q}}^{\orbify{P}}\cong N_{\Gamma_\orbify{P}^\orbify{O}}(\Delta)/C_{\Gamma_\orbify{P}^\orbify{O}}(\Delta)$. For convenience, we will denote $\Gamma_\orbify{P}^\orbify{O}$ by $\Gamma_{\orbify{P}}\cong A_4$. Now, since $\Delta\cong\ZZ_3$, the order of the normalizer $N_{\Gamma_\orbify{P}}(\Delta)$ must be 3, 6, or 12. Since there are no subgroups of order 6 in $A_4$, and $\ZZ_3$ is not normal in $A_4$, it must be the case that $N_{\Gamma_\orbify{P}}(\Delta)\cong\ZZ_3$. Thus, we also conclude that $C_{\Gamma_\orbify{P}}(\Delta)\cong\ZZ_3$, and finally that
$$\Gamma_{\orbify{Q}}^{\orbify{P}}\cong N_{\Gamma_\orbify{P}}(\Delta)/C_{\Gamma_\orbify{P}}(\Delta)\cong\{1\}.$$
This means that in the canonical orbifold substructure that $\orbify{Q}$ inherits from $\orbify{P}$, the origin has trivial isotropy.

So, we have shown that the inherited orbifold substructure, $\orbify{Q}^\orbify{O}$, that $\orbify{Q}$ inherits from $\orbify{O}$ is different than, $\orbify{Q}^\orbify{P}$, the one that $\orbify{Q}$ inherits from $\orbify{P}$. It is useful to note, however, that it is possible to recover $\orbify{Q}^\orbify{P}$ as a suborbifold of $\orbify{O}$, by choosing $\Lambda=A_3\cong\ZZ_3$ as opposed to the entire stabilizer group $N_\Gamma(\Delta)\cong\ZZ_2\ltimes\ZZ_3$. That is, $\orbify{Q}^\orbify{P}$ is, in fact, a suborbifold of $\orbify{O}$ with a non-canonical inherited substructure!

\section{Inherited Canonical Orbifold Substructures and Saturation}\label{CanVsSatSection}
In this section we show that the suborbifold $\orbify{P}\subset\orbify{O}$ constructed in section~\ref{MainCounterExample} is not a saturated suborbifold even though $\orbify{P}$ has the canonical orbifold substructure inherited from $\orbify{O}$. Recall, $\tilde W=\RR[\Delta]=\RR[\ZZ_3]=\text{span}\{1,\gamma,\gamma^2\}\subset\tilde V=\RR[\Beta]=\RR[A_4]$, and let $\tilde w=c_11+c_\gamma\gamma+c_{\gamma^2}\gamma^2$. Let $\beta\in\Beta\cong A_4$. Then 
$$\beta\cdot \tilde w=c_11+c_\gamma\gamma^\beta+c_{\gamma^2}(\gamma^\beta)^2$$
where $\gamma^\beta=\beta\gamma\beta^{-1}$. Then, either $\gamma^\beta\notin\text{span}\{1,\gamma,\gamma^2\}$ or $\gamma^\beta=\gamma$. The second case follows because we know the normalizer $N_\Beta(\Delta)\cong\ZZ_3$. Thus, the orbit $(A_4\cdot\tilde w)\cap\tilde W=\tilde w$. Now let 
$$\alpha\in N_\Gamma(\Delta)=N_{A_5}(\ZZ_3)\cong\ZZ_2\ltimes\ZZ_3= \langle \alpha,\gamma\mid \alpha^2=\gamma^3=1, \gamma^\alpha=\gamma^{-1}\rangle.$$
Then $\alpha\notin A_4$ and
$$\alpha\cdot\tilde w=c_11+c_\gamma\gamma^\alpha+c_{\gamma^2}(\gamma^\alpha)^2=c_11+c_\gamma\gamma^2+c_{\gamma^2}\gamma\in\tilde W.$$
This implies that $(A_5\cdot\tilde w)\cap\tilde V\ne A_4\cdot\tilde w$, and thus $\orbify{P}$ is not saturated in $\orbify{O}$.

\end{document}